\documentclass[12pt,english]{extarticle}
\usepackage[T1]{fontenc}
\usepackage[latin9]{inputenc}
\usepackage{color}
\usepackage{babel}
\usepackage{amsmath}
\usepackage{amsthm}
\usepackage{amssymb}
\usepackage{geometry}
\usepackage{microtype}
\usepackage[pdfusetitle,
 bookmarks=true,bookmarksnumbered=true,bookmarksopen=false,
 breaklinks=false,pdfborder={0 0 1},backref=false,colorlinks=true]
 {hyperref}
\hypersetup{
 linkcolor=blue, urlcolor=blue, citecolor=blue, pdfstartview={FitH}, hyperfootnotes=false, unicode=true}

\makeatletter
\numberwithin{equation}{section}
\numberwithin{figure}{section}
\theoremstyle{plain}
\newtheorem{thm}{\protect\theoremname}
\theoremstyle{definition}
\newtheorem{defn}[thm]{\protect\definitionname}
\theoremstyle{remark}
\newtheorem{rem}[thm]{\protect\remarkname}

\usepackage{lettrine} 
\let\myFoot\footnote
\renewcommand{\footnote}[1]{\myFoot{#1\vspace{3mm}}}
\usepackage{amsmath}
\usepackage{amssymb}
\usepackage{xcolor}
\usepackage{pdfrender}
\usepackage{mlmodern}
\usepackage[T1]{fontenc}

\usepackage{amsmath}

\usepackage{tikz}
\usepackage{circuitikz}

\makeatother

\providecommand{\definitionname}{Definition}
\providecommand{\remarkname}{Remark}
\providecommand{\theoremname}{Theorem}

\begin{document}
\title{A note on the planar Skorokhod embedding problem}
\author{Maher Boudabra \thanks{King Fahd University of Petroleum and Minerals, KSA. https://orcid.org/0000-0001-5043-0058}}
\maketitle
\begin{abstract}
The planar Skorokhod embedding problem was first proposed and solved
by R. Gross in $2019$ \cite{gross2019}. Gross worked with probability
distributions having finite second moment. In \cite{boudabra2019remarks,Boudabra2020},
the solutions extended to all distributions with a finite $p^{th}$
moment for $p>1$. The case $p=1$ remained uncovered since then.
In this note we show that the planar Skorokhod embedding problem is
solvable for $p=1$ when the Hilbert transform of its quantile function
is integrable, effectively closing this line of investigation.
\end{abstract}
MSC2020 subject classifications: 30C35; 60D05; 51M25

\section{Introduction and main results}

In $2019$, R. Gross posed the following problem \cite{gross2019}:
Given a non-degenerate distribution $\mu$ with zero mean and finite
second moment, find a simply connected domain $U$ such that if $(Z_{t})_{t\geq0}$
is a planar Brownian motion started at $0$  and $\tau_{U}$ is its
exit time from $U$, then $\Re(Z_{\tau_{U}})$ has the law $\mu$. 

This is called the planar (or conformal) Skorokhod embedding problem
(PSEP). Gross gave a constructive solution using the principle of
conformal invariance of planar Brownian motion \cite{revuz2013continuous}.
The domain he constructed enjoys the property that its underlying
exit time has a finite average. In \cite{boudabra2019remarks}, Boudabra
and Markowsky showed that Gross' approach remains valid if the distribution
$\mu$ has a finite $p^{th}$ moment for some $p>1$. They reformulated
the problem by requiring the finiteness of $\mathbb{E}(\tau_{U}^{\frac{p}{2}})$.
We adopt this formulation in this note. Later on, in a subsequent
work \cite{Boudabra2020}, the same authors provided a second solution
to the PSEP with a different geometry from Gross' domain. Both solutions
start by crafting a certain univalent function acting on the open
unit disc $\mathbb{D}$. Then the authors show that the range $U$
of such a function fits the statistical requirement $\Re(Z_{\tau_{U}})\sim\mu$.
Our main result gives a sufficient condition for the existence of
a $p=1$ embedding in terms of the Hilbert transform of the quantile
function; see Theorem \ref{thm:main} below. Through this work, the
following notations are adopted: 
\begin{itemize}
\item $\mathbb{D}$: the open unit disc in the plane.
\item $\mathbb{T}$: the unit circle in the plane.
\item $L_{2\pi}^{p}$: the space of real valued $2\pi$-periodic functions
$f$ such that
\[
\Vert f\Vert_{p}=\left(\frac{1}{2\pi}\int_{0}^{2\pi}\vert f(t)\vert^{p}dt\right)^{\frac{1}{p}}
\]
 is finite ($p>0$). We write $L^{p}(\mathbb{T})$ for the space of
complex-valued functions on $\mathbb{T}$. 
\item Standard planar Brownian motion: planar Brownian motion started at
$0$.
\end{itemize}
The ingredients of the proofs of the results of our work are Hardy
spaces and Hilbert transform. We recall their definitions for the
reader's convenience. 
\begin{defn}
Let $f$ be an analytic function in the unit disc. The $p^{th}$ Hardy
norm of $f$ is defined by 
\begin{equation}
\Vert f\Vert_{\boldsymbol{H}^{p}}:=\sup_{0\leq r<1}\left\{ \frac{1}{2\pi}\int_{0}^{2\pi}|f(re^{\theta i})|^{p}d\theta\right\} ^{\frac{1}{p}}.\label{hardy}
\end{equation}
The set of analytic functions on $\mathbb{D}$ with finite $\boldsymbol{H}^{p}$
norm is denoted by $\boldsymbol{H}^{p}(\mathbb{D})$ and is called
the Hardy space (of index $p$).
\end{defn}

\begin{defn}
The Hilbert transform of a $2\pi$-periodic function $f$ is defined
by 

\[
\mathcal{H}\{f\}(x):=\frac{1}{2\pi}P.V\left\{ \int_{-\pi}^{\pi}f(x-t)\cot({\textstyle \frac{t}{2}})dt\right\} =\frac{1}{2\pi}\lim_{\eta\rightarrow0}\int_{\eta\leq|t|\leq\pi}f(x-t)\cot({\textstyle \frac{t}{2}})dt
\]
where $P.V$ denotes the Cauchy principal value. Note that the Hilbert
transform is an anti-involution on zero-mean functions, i.e. $\mathcal{H}^{2}=-id$.
The Hilbert transform exists a.e. for any function $f\in L^{p}(\mathbb{T})$
with $p\geq1$. In addition, if $p>1$ then 
\begin{equation}
\Vert H\{f\}\Vert_{p}\leq\lambda_{p}\Vert f\Vert_{p}\label{strong inequality-1}
\end{equation}
for some positive constant $\lambda_{p}$. See \cite{duren2000theory}.
\end{defn}

The main results of the work are the following.
\begin{thm}
\label{boundary distribution} Let $f$ be a univalent function on
$\mathbb{D}$ with $f(0)=0$ and set $U=f(\mathbb{D})$. Let $(Z_{t})_{t\geq0}$
be standard planar Brownian motion  and $\tau_{U}$ its exit time
from $U$. Then $Z_{\tau_{U}}$ and $f(\xi)$ share the same distribution
where $\xi$ is any random variable uniformly distributed on the unit
circle and $f(\xi)$ denotes the radial boundary value of $f$ at
$\xi$ (defined a.e.).
\end{thm}

\begin{thm}
\label{thm:main} Assume $\mu$ is a non-degenerate probability measure
with zero mean and let $q$ be its quantile function. Then there exists
a simply connected domain $U\subset\mathbb{C}$ containing the origin
such that, for a standard planar Brownian motion $(Z_{t})_{t\ge0}$
 and $\tau_{U}:=\inf\{t\ge0:Z_{t}\notin U\}$, we have 
\[
\Re(Z_{\tau_{U}})\sim\mu\,\,\,\text{and}\,\,\,\mathbb{E}(\tau_{U}^{\frac{1}{2}})<+\infty
\]
whenever $\mathcal{H}\{q(\frac{\vert\cdot\vert}{\pi})\}\in L_{2\pi}^{1}$
where $q$ is the quantile function of $\mu$. 
\end{thm}

\section{Proofs}

Before proving our main results, we recall two celebrated results
that we need in the proof. 
\begin{thm}
\label{thm:p<1/2}\cite[Theorem 3.16]{duren2001univalent} Any univalent
function acting on the unit disc belongs to Hardy spaces $\boldsymbol{H}_{p}(\mathbb{D})$
for all $p<\frac{1}{2}$.
\end{thm}

\begin{thm}
\cite{Rudin2001} \label{hardy radial norm} If $f\in\boldsymbol{H}_{p}(\mathbb{D})$
for some positive $p$ then 
\begin{enumerate}
\item the radial limit of $f$, say $f^{\ast}(\xi)$, exists a.e on $\mathbb{T}$
and $f^{*}\in L^{p}(\mathbb{T})$.
\item $\underset{r\rightarrow1^{-}}{\lim}\Vert f(r\cdot)-f^{*}\Vert_{p}=0$.
\item $\Vert f^{*}\Vert_{p}$ equals the $p^{th}$ Hardy norm of $f$.
\end{enumerate}
\end{thm}

From now on, we do not distinguish notationally between a Hardy function
and its boundary values; $f(\xi)$ will denote the radial limit on
$\mathbb{T}$.

\subsection{Proof of theorem \ref{boundary distribution}}

The proof is adapted from Lemma 3.2 in \cite{boudabra2019remarks}.
The main difference is that in \cite{boudabra2019remarks}, the authors
required that $f$ belongs to some Hardy space $\boldsymbol{H}_{p}(\mathbb{D})$
for some $p>1$, while this constraint is dropped here. 

Let $(W_{t})_{t}$ be a standard planar Brownian motion and killed
on exiting $\mathbb{D}$. By the conformal invariance principle, the
planar process $\left(Z_{t}=f(W_{\sigma^{-1}(t)})\right)_{t\geq0}$
is a standard planar Brownian motion, where $\sigma$ is the time
change given by 
\[
\sigma(t)=\int_{0}^{t}\vert f'(W_{s})\vert^{2}ds,\,\,t<\tau_{\mathbb{D}}.
\]
In particular 
\[
Z_{\sigma(\tau_{n})}=f(W_{\tau_{n}})
\]
where $\tau_{n}$ denotes the exit time of $W_{t}$ from $\mathbb{D}_{n}:=(1-\frac{1}{n})\mathbb{D}$.
In fact $\sigma(\tau_{n})$ is the exit time of $(Z_{t})_{t\geq0}$
from $f(\mathbb{D}_{n}).$ Since $\mathbb{D}_{n}\uparrow\mathbb{D}$
as $n\rightarrow+\infty$, we have $\tau_{n}\uparrow\tau_{\mathbb{D}}$
almost surely. By continuity of the time change, $\sigma(\tau_{n})$
converges to $\sigma(\tau_{\mathbb{D}})$ a.s. The quantity $\sigma(\tau_{\mathbb{D}})$
is exactly the exit time of $(Z_{t})_{t\geq0}$ from $U$, which we
denote by $\tau_{U}$. As the Brownian path is continuous a.s., $Z_{\sigma(\tau_{n})}$
converges to $Z_{\tau_{U}}$ a.s. On an auxiliary probability space
$\Omega=(0,2\pi)$ equipped with the normalized Lebesgue measure,
define $\xi_{n}(\omega)=(1-\frac{1}{n})e^{\omega i}$ and $\xi(\omega)=e^{\omega i}$
where $\omega\in\Omega$. As planar Brownian motion is isotropic,
the two random variables $W_{\tau_{n}}$ and $\xi_{n}(\omega)=(1-\frac{1}{n})e^{\omega i}$
share the same distribution. In particular $f(\xi_{n})$ and $Z_{\sigma(\tau_{n})}$
have the same law as well. Now fix $p\in(0,\frac{1}{2})$. By Theorem
\ref{thm:p<1/2}, $f$ belongs to $\boldsymbol{H}_{p}(\mathbb{D})$
and hence it extends radially to the boundary of $\mathbb{D}$. We
have
\[
\mathbb{E}(\vert f(\xi_{n})-f(\xi)\vert^{p})=\frac{1}{2\pi}\int_{0}^{2\pi}|f((1-{\textstyle \frac{1}{n}})e^{\omega i})-f(e^{\omega i})|^{p}d\omega=\Vert f((1-{\textstyle \frac{1}{n}})\cdot)-f\Vert_{p}^{p}.
\]
By the second point of Theorem \ref{hardy radial norm}, $\mathbb{E}(\vert f(\xi_{n})-f(\xi)\vert^{p})$
converges to zero. Hence $f(\xi_{n})$ converges to $f(\xi)$ in distribution.
As $f(\xi_{n})$ and $Z_{\sigma(\tau_{n})}$ are equal in distribution,
$Z_{\sigma(\tau_{n})}$ converges to $f(\xi)$ in distribution as
well. On the other hand 
\[
\lim_{n}Z_{\sigma(\tau_{n})}=Z_{\tau_{U}}\,\,a.s.
\]
 Therefore $Z_{\tau_{U}}$ and $f(\xi)$ have the same distribution
by uniqueness of weak limits. This ends the proof. 

\subsection{Proof of theorem \ref{thm:main}}

Recall that the quantile function $q$ of $\mu$ is simply the pseudo
inverse of the cumulative distribution function $F$ of $\mu$. More
precisely 
\[
q(u)=\inf\{x\mid F(x)\geq u\}.
\]
A crucial property of $q$ is that when fed with a uniformly distributed
random variable $X$ in $(0,1)$ then $q(X)$ samples as $\mu$. Set
$\varphi(\theta)=q(\frac{\vert\theta\vert}{\pi})$ with $\theta\in(-\pi,\pi)$.\\
\\
The assumption $H\{\varphi\}\in L_{2\pi}^{1}$ implies that the analytic
map 
\[
G(z)=\frac{1}{2\pi}\int_{0}^{2\pi}\frac{e^{ti}+z}{e^{ti}-z}\varphi(t)dt\in\boldsymbol{H}^{1}(\mathbb{D}).
\]
The map $G$ is the standard Cauchy-Poisson integral (See \cite[Chapter 3]{duren2000theory}).
In particular 
\begin{equation}
\Re(G(e^{\theta i}))=\varphi(\theta),\,\,\Im(G(e^{\theta i}))=H\{\varphi\}\,\,\,a.e.\label{real of boundary}
\end{equation}
Furthermore, the map $G(z)$ is univalent by Proposition $2.2$ in
\cite{gross2019}. In particular, $U:=G(\mathbb{D})$ is simply connected
and contains the origin. Run a planar Brownian motion $\left(Z_{t}\right)_{t}$
inside $U$ and set $\tau_{U}$ the exit time of $Z_{t}$ from $U$.
By \ref{real of boundary}, the boundary value of $G$ has real part
$\varphi$, which encodes the quantile of $\mu$. Therefore, by Theorem
\ref{boundary distribution}, 
\[
\Re(Z_{\tau_{U}})\sim\mu.
\]
Since $G\in\boldsymbol{H}^{1}(\mathbb{D})$, Burkholder's characterization
of the finiteness of the moments of $\tau_{U}$ \cite{burkholder1977exit},
yields $\mathbf{E}(\tau_{U}^{\frac{1}{2}})<+\infty$. Therefore, the
PSEP is solvable. \\
\\

\section{Comments}

\subsubsection*{Comment 1}

Theorem \ref{boundary distribution} might seem obvious as one may
think that
\[
\lim_{n}f(W_{\tau_{n}})=f(W_{\tau_{\mathbb{D}}})
\]
is immediate. But this is not obvious for the simple reason that $W_{\tau_{n}}$
does not approach the boundary of the unit disc radially! Theorem
\ref{thm:main} naturally raises the converse question: if the PSEP
is solvable for a given law $\mu$, must the Hilbert transform $\mathcal{H}\{q(\frac{\vert\cdot\vert}{\pi})\}\in L_{2\pi}^{1}$? We conjecture that the answer is yes. 

\subsubsection*{Comment 2}

Theorem \ref{thm:main} implies the existence of solutions for the
PSEP when $\mu$ has a finite $p^{th}$ moment for some $p>1$. This
follows from \ref{strong inequality-1}. \\
\\
A natural question is how to find functions with integrable Hilbert
transforms in $L_{2\pi}^{1}$. A naive way to do so is to compute
the Hilbert transform and then check its integrability. However, this
is hard in practice. The following result, due to Zygmund, is handy.
\begin{thm}
\cite{duren2000theory} Let $u$ be a harmonic function in the unit
disc such that
\begin{equation}
\sup_{0\leq r<1}\int_{0}^{2\pi}\vert u(re^{\theta i})\vert\log^{+}(\vert u(re^{\theta i})\vert)d\theta<+\infty\label{zygmund}
\end{equation}
then 
\[
\sup_{0\leq r<1}\int_{0}^{2\pi}\vert H\{u\}(re^{\theta i})\vert d\theta<+\infty.
\]
\end{thm}

We invite the reader to have a look at the valuable discussion of
Zygmund's theorem provided in \cite[Chapter 4]{duren2000theory}.
In particular, the condition \ref{zygmund} is, in some sense, the
weakest one can impose on $u$. Note that the $\log^{+}(x)$ can be
replaced by $\log(1+\vert x\vert)$ \cite{garnett2006bounded}. A
classical result in Harmonic Analysis says that if $g$ is a convex
function then $g(u)$ is subharmonic. In particular, the integral
in \ref{zygmund} is non-decreasing in $r$ \cite{Rudin2001}. So,
in terms of boundary function, the condition \ref{zygmund} can be
replaced by 
\[
\int_{0}^{2\pi}\vert u(e^{\theta i})\vert\log^{+}(\vert u(e^{\theta i})\vert)d\theta<+\infty,
\]
which was originally formulated. Based on this fact, the planar Skorokhod
embedding problem becomes solvable once $\vert\varphi\vert\log^{+}(\vert\varphi\vert)$
belongs to $L_{2\pi}^{1}$. 

Now we provide an example where the distribution $\mu$ has only a
first moment. Consider the distribution of the random variable 
\[
\xi=\frac{5}{(\ln(10)+1)^{2}}-\frac{10}{X(1-\ln(\frac{X}{10}))^{3}}
\]
 with $X$ uniformly distributed on $(0,1)$. The underlying quantile
function is 
\[
q(u)=\frac{5}{(\ln(10)+1)^{2}}-\frac{10}{u(1-\ln(\frac{u}{10}))^{3}},\,\,u\in(0,1).
\]
Using Bertrand's criterion for integrals, $q$ belongs to $L^{p}$
for $p=1$ but not for any higher $p$. Near zero, we have the approximation
\begin{equation}
\frac{\log^{+}(\frac{1}{u(1-\ln(u))^{3}})}{u(1-\ln(u))^{3}}=-\frac{\log(u(1-\ln(u))^{3})}{u(1-\ln(u))^{3}}\approx-\frac{\log(u)}{u\vert\ln(u)\vert^{3}}-3\frac{\log(\vert\ln(u)\vert)}{u\vert\ln(u)\vert^{3}}.\label{bertrand}
\end{equation}
Thus, the integral of $\vert q\vert\log^{+}(\vert q\vert)$ is finite.
Therefore the Hilbert transform of $\varphi$ is in $L_{2\pi}^{1}$
(recall $\varphi(\theta)=q(\frac{\vert\theta\vert}{\pi})$). Consequently,
the Skorokhod embedding problem is solvable for the considered $\mu$. 
\begin{rem}
Gross' construction of the map $G(z)$ is based on the calculation
of the Fourier coefficients of $\varphi$. Now one can see that $G(z)$
can be simply expressed as
\[
G(z)=\frac{1}{2\pi}\int_{0}^{2\pi}\frac{e^{ti}+z}{e^{ti}-z}\varphi(t)dt.
\]
Although this paper is primarily addressing the method introduced
in \cite{gross2019}, the approach is directly applicable, without
modification, to the method described in \cite{Boudabra2020}. 
\end{rem}

\subsubsection*{Comment 3}

Mathematicians tend to shy away from $L^{p}$ spaces when $p<1$ due
to their unconventional properties, such as the failure of the triangle
inequality, or, in relation to our context, the inability to construct
Poisson integrals. However, we present an example that shows something
might be done when $p<1$. Using the formulas in \cite{boudabra2019remarks,boudabra2021some},
the c.d.f of the stopped planar Brownian motion upon leaving the Koebe
domain $\mathcal{K}=\mathbb{C}-(-\infty,-\frac{1}{4}]$ is given by
\[
F(w)=1-{\textstyle \frac{2}{\pi}}\arctan(\sqrt{-1-4w})
\]
for $w\in(-\infty,-\frac{1}{4}]$. That is, we get 
\[
\varphi(t)=F^{-1}({\textstyle \frac{\vert t\vert}{\pi}})=-\frac{1}{4\sin(\frac{t}{2})^{2}}.
\]
Interestingly, the function $\varphi(t)$ coincides with the real
part of the radial limit of the Koebe function $k(z)=\frac{z}{(1-z)^{2}}$.
In \cite{gardiner2006radial}, the author gathered some results about
the existence of entire functions and harmonic functions with a predefined
radial boundary function. So a natural question pops up and we think
it is worth being addressed: How to construct the univalent function
whose real part of its radial limit matches $\varphi$ when the distribution
$\mu$ has exclusively a finite $p^{th}$ moment for some $p<1$?
Another related issue will be how to reconcile this with the fact
that every univalent function belongs to $\boldsymbol{H}^{p}(\mathbb{D})$
for all $p<\frac{1}{2}$? 

\subsubsection*{Comment 4}

It is worth noting that Theorem \ref{boundary distribution} is still
valid if $f$ is a proper map, or more generally a $B$-proper map.
The notion of $B$-properness was introduced in \cite{markowsky2015exit}.
Both types of functions ensure that if a sequence $(z_{n})_{n}$ in
the input domain is approaching its boundary then the image sequence
$(f(z_{n}))_{n}$ is also approaching the boundary of the range of
$f$. 

\subsubsection*{Acknowledgment}

I would like to thank Professor Zihua Guo from Monash University for
his valuable comments. 

\bibliographystyle{plain}
\bibliography{MaherLibrary}

@Book{revuz2013continuous,
  author    = {D. Revuz and M. Yor},
  publisher = {Springer Science \& Business Media},
  title     = {Continuous martingales and Brownian motion},
  year      = {2013},
  volume    = {293},
}

@Article{burkholder1977exit,
  author    = {Burkholder, D.},
  title     = {Exit times of {B}rownian motion, harmonic majorization, and {H}ardy spaces},
  journal   = {Advances in {M}athematics},
  year      = {1977},
  volume    = {26},
  number    = {2},
  pages     = {182--205},
  publisher = {Academic Press},
}

@Book{duren2000theory,
  author    = {P. L Duren},
  publisher = {Courier Corporation},
  title     = {Theory of $H^p$ spaces},
  year      = {2000},
}

@Article{gross2019,
  author  = {Gross, R.},
  title   = {A conformal {S}korokhod embedding},
  journal = {Electronic Communications in Probability},
  year    = {2019},
}

@Book{Rudin2001,
  author    = {Rudin, W.},
  publisher = {McGraw-Hill Education},
  title     = {Real and complex analysis (3rd ed.)},
  year      = {2001},
}

@Article{boudabra2019remarks,
  author  = {Boudabra, M. and Markowsky, G.},
  journal = {Electronic Communications in Probability},
  title   = {Remarks on {G}ross' technique for obtaining a conformal {S}korokhod embedding of planar {B}rownian motion},
  year    = {2020},
}

@Article{Boudabra2020,
  author  = {M. Boudabra and G. Markowsky},
  journal = {Journal of Mathematical Analysis and Applications},
  title   = {A new solution to the conformal {S}korokhod embedding problem and applications to the {D}irichlet eigenvalue problem},
  year    = {2020},
  issn    = {0022-247X},
  number  = {2},
  pages   = {124351},
  volume  = {491},
}

@Article{markowsky2015exit,
  author    = {G. Markowsky},
  journal   = {Journal of Mathematical Analysis and Applications},
  title     = {The exit time of planar {B}rownian motion and the {P}hragm{\'e}n--{L}indel{\"o}f principle},
  year      = {2015},
  number    = {1},
  pages     = {638--645},
  volume    = {422},
  publisher = {Elsevier},
}

@Book{duren2001univalent,
  author    = {P. L Duren},
  publisher = {Springer Science \& Business Media},
  title     = {Univalent functions},
  year      = {2001},
  volume    = {259},
}

@PhdThesis{boudabra2021some,
  author = {M. Boudabra},
  school = {Monash University},
  title  = {On some problems related to exit times of planar Brownian motion},
  year   = {2021},
}

@InCollection{gardiner2006radial,
  author    = {S J. Gardiner},
  booktitle = {Potential Theory in Matsue},
  publisher = {Mathematical Society of Japan},
  title     = {Radial limits of harmonic functions},
  year      = {2006},
  pages     = {43--52},
  volume    = {44},
}

@Book{garnett2006bounded,
  author    = {J. Garnett},
  publisher = {Springer Science \& Business Media},
  title     = {Bounded analytic functions},
  year      = {2006},
  volume    = {236},
}

\end{document}